\documentclass[12pt]{amsart}
\usepackage{times,a4wide}
\usepackage{amssymb}

\usepackage{color}
\usepackage{epsfig}

\newtheorem*{theorem*}{Theorem}

\newtheorem{remark}{Remark}
\newtheorem*{remark*}{Remark}

\newtheorem{theorem}{Theorem}
\newtheorem{proposition}{Proposition}
\newtheorem*{proposition*}{Proposition}
\newtheorem{corollary}{Corollary}
\newtheorem*{corollary*}{Corollary}
\newtheorem{lemma}{Lemma}

\newtheorem*{mlemma*}{Main Lemma}
\newtheorem*{lemma*}{Lemma}

\newtheorem*{claim*}{Claim}
\theoremstyle{remark}
\newtheorem*{example}{Example}

\newcommand{\nc}{\newcommand}
\newcommand{\rnc}{\renewcommand}

\newcommand{\C}{{\mathbb C}}
\newcommand{\N}{{\mathbb N}}
\newcommand{\D}{{\mathbb D}}

\newcommand{\T}{{\mathbb T}}
\newcommand{\Z}{{\mathbb Z}}

\nc{\ol}{\overline}
\nc{\bea}{\begin{eqnarray}}
\nc{\eea}{\end{eqnarray}}
\nc{\beqa}{\begin{eqnarray*}} \nc{\bega}{\begin{gather*}}
\nc{\Sd}{\stackrel{\cdot}{+}} \nc{\eeqa}{\end{eqnarray*}}
\nc{\eega}{\end{gather*}} \nc{\nn}{\nonumber}
\nc{\lmt}{\longmapsto} \nc{\eps}{\varepsilon}
\nc{\Llra}{\Longleftrightarrow} \nc{\LRA}{\Leftrightarrow}
\nc{\lra}{\longrightarrow} \nc{\Lra}{\Longrightarrow}
\nc{\Lla}{\Longleftarrow} \nc{\lmto}{\longmapsto}
\nc{\vt}{\vartheta} \nc{\vs}{\psi} \nc{\sn}{\sigma_n}
\nc{\Hi}{H^{\infty}} \nc{\Li}{L^{\infty}} \nc{\vp}{\varphi}
\nc{\Tp}{T_{\vp}} \nc{\Tpb}{T_{\overline{\vp}}}
\nc{\Tpn}{T_{\vp_0}} \nc{\Tpnb}{T_{\overline{\vp}_0}}
\nc{\gp}{g^{2/p}} \nc{\Kpw}{K^p_I(|g|^2)} \nc{\Kqw}{K^q_I(|g|^2)}
\nc{\Kp}{K^p_I} \nc{\kep}{\operatorname{ker} T_{\vp}}
\nc{\Lpw}{L^p(|g|^2)} \nc{\Lppw}{L^{p'}(|g|^2)}
\nc{\Lqw}{L^q(|g|^2)} \nc{\Ktw}{K^2_I(|g|^2)} \nc{\Kt}{K^2_I}
\nc{\Ltw}{L^2(|g|^2)} \nc{\kl}{k_{\lambda}}
\nc{\klb}{k^b_{\lambda}} \nc{\klIb}{k^{Ib}_{\lambda}}
\nc{\klI}{k_{\lambda}^I} \nc{\klIt}{k_{\lambda}^{\tilde{I}/I}}

\nc{\Rg}{\operatorname{Rg}} \nc{\essinf}{\operatorname{essinf}}
\rnc{\Re}{\operatorname{Re}} 
\nc{\Lin}{\operatorname{Lin}} \nc{\const}{\operatorname{const}}
\nc{\dist}{\operatorname{dist}} \nc{\dpst}{\displaystyle}
\nc{\Hol}{\operatorname{Hol}}
\nc{\clos}{\operatorname{clos}}

\nc{\cl}{c_{\lambda}}
\nc{\dl}{d_{\lambda}}
\nc{\el}{e_{\lambda}}

\nc{\NM}[1]{\|#1\|_{\mathcal{M}(A)}}
\nc{\NH}[1]{\|#1\|_{\mathcal{H}(A)}}
\nc{\NMp}[2]{\|#1\|_{\mathcal{M}(#2)}}
\nc{\NHp}[2]{\|#1\|_{\mathcal{H}(#2)}}

\nc{\Mph}{{\mathcal{M}(\vp)}}
\nc{\Hph}{{\mathcal{H}(\vp)}}

\nc{\Mp}{\mathcal{M}}
\nc{\Hp}{\mathcal{H}}

\renewcommand{\Re}{\operatorname{Re}}
\renewcommand{\Im}{\operatorname{Im}}

\title[Analytic continuation and embeddings]{Analytic continuation and
  embeddings in weighted backward shift invariant subspaces}

\author{Andreas Hartmann}

\keywords{analytic continuation, backward shift, embedding,
invariant subspace, Muckenhoupt
condition, spectrum, Toeplitz operators}

\subjclass[2000]{30B40,30D55,47A15,47B37}

\begin{document}

\begin{abstract}
By a famous result, 
functions in backward shift invariant
subspaces in Hardy spaces are characterized by the fact that
they admit a pseudocontinuation a.e.\ on $\T$. More can be
said if the spectrum of the associated inner function has
holes on $\T$. Then the functions of the invariant subspaces
even extend analytically through these holes. We will discuss
the situation in weighted backward shift invariant subspaces.
The results on analytic continuation will be applied  
to consider some embeddings of weighted invariant subspaces into 
their unweighted companions.
Such weighted versions of invariant subspaces appear 
naturally in the context of Toeplitz operators.
A connection between the spectrum of the inner function and 
the approximate point spectrum of the backward shift in the
weighted situation is established
in the spirit of results by Aleman, Richter and Ross. 
\end{abstract}

\maketitle

\section{Introduction}\label{S1}

Backward shift invariant subspaces have shown to be of great interest in many 
domains in complex analysis and Operator Theory. In 
$H^2$, the classical Hardy space of holomorphic functions on
the unit disk $\D$ satisfying
\beqa
 \|f\|_2^2:=\sup_{0<r<1}\frac{1}{2\pi}\int_{-\pi}^{\pi}
 |f(re^{it})|^2dt<\infty,
\eeqa
they are given by $H^2\ominus IH^2$, where $I$ is an inner
function, that is a bounded analytic function in $\D$ the
boundary values of which are in modulus equal to 1 a.e.\ on $\T$.
Another way of writing the model spaces is
\beqa
 K^2_I=H^2\cap I\overline{H^2_0},
\eeqa
where $H^2_0=zH^2$ is the subspace of functions in $H^2$
vanishing in 0. The bar sign means complex conjugation here. 
This second writing $K^2_I=H^2\cap I\overline{H^2_0}$
does not appeal to the Hilbert space structure and 
thus generalizes
to $H^p$ (which is defined
as $H^2$ but replacing the integration power $2$ by $p\in (0,\infty)$;
it should be noted that for $p\in (0,1)$ the expression $\|f\|_p^p$
defines a metric;
for $p=\infty$, $\Hi$ is the Banach space of bounded analytic
functions on $\D$ with obvious norm).
When $p=2$, then
these spaces are also called model spaces.
Model spaces have attracted a lot of attention
of course 
in operator theory, initially in the function model of Nagy 
and Foias, but then also in perturbation theory with Clark's
seminal paper on rank one perturbations of the compressed shift
on $K^2_I$. 
As a result of Clark the Cauchy transform
allows to identify $K^2_I$ with $L^2(\sigma_{\alpha})$
where $\sigma_{\alpha}$ is a so-called Clark measure that 
one can deduce from $I$. Clark's motivation was in fact
to consider completeness problems in model spaces $K^2_I$.
In a series of papers, Aleksandrov and Poltoratski
were interested especially in the behaviour of the Cauchy transform
when $p\neq 2$. 

Another interest in backward shift invariant subspaces concerns
embedding questions, especially when $K^p_I$ embeds into
some $L^p(\mu)$. Those questions were investigated for
instance by Aleksandrov, Treil, Volberg and many others
(see for instance \cite{TV1} for some results). Here we
will in fact be interested in the different situation when
the weight is not on $L^p$ but on $K^p_I$.

A very important result in connection with $K^p_I$-space
is that of Douglas, Shapiro and Shields (\cite{DSS},
see also \cite[Theorem 1.0.5]{CR}). 
They have in fact characterized
$K^p_I$ as the space of functions in $H^p$ that admit a
so-called pseudocontinuation. Recall that a function holomorphic
in $\D_e:=\hat{\C}\setminus{\clos\D}$ --- we will use $\clos E$ to
designate the closure of a set $E$ in order to preserve the bar-sign
for complex conjugation --- is a pseudocontinuation of a 
function $f$ meromorphic in $\D$ if $\psi$ vanishes at
$\infty$ and the outer nontangential limits of $\psi$ on $\T$
coincide with the inner nontangential limits of $f$ on $\T$
in almost every point of $\T$. Note that $f\in K^2_I=H^2\cap
I\overline{H^2_0}$ implies that $f=I\overline{\psi}$ with
$\psi\in H^2_0$. Then the meromorphic function $f/I$
equals $\overline{\psi}$ a.e.\ $\T$, and writing 
$\psi(z)=\sum_{n\ge 1}b_nz^n$,
it is clear that $\tilde{\psi}(z):=\sum_{n\ge 1}\overline{b_n}/z^n$
is a holomorphic function in $\D_e$, vanishing at $\infty$, and
being equal to $f/I$ almost everywhere on $\T$
(in fact, $\tilde{\psi}\in H^2(\D_e)$). 

Note that there are functions analytic on $\C$ that do not
admit a pseudocontinuation. An example of such a function is 
$f(z)=e^z$ which has an essential singularity at infinity.

On the other hand, there are of course pseudocontinuations that are
not analytic continuations. A 
result by Moeller \cite{Mo} states that outside the spectrum of $I$,
$\sigma(I)=\{\lambda\in \clos\D:\liminf_{z\to\lambda} I(z)=0\}$,
which is  a closed set, every function $f\in K^2_I$
extends analytically through the circle. It is not difficult
to construct inner functions $I$ for which $\sigma(I)\cap \T
=\T$. Take for instance for $I$ the Blaschke product
associated with the sequence $\Lambda=\{(1-1/n^2)e^{in}\}_n$, the
zeros of which accumulate at every point on $\T$.
\\

The problem we are interested in is the case of a weighted
backward shift invariant subspace. Let $I$ be any inner function, and
$g$ an outer function in $H^p$, $1<p<\infty$. Set
\beqa
 K^p_I(|g|^p)=H^p(|g|^p)\cap I\overline{H^p_0(|g|^p)}.
\eeqa
Here 
\beqa
 H^p(|g|^p)=\{f\in \Hol(\D):\|f\|^p_{|g|^p}&:=&\sup_{0<r<1}\frac{1}{2\pi}
 \int_{-\pi}^{\pi} |f(re^{it})|^p|g(re^{it})|^pdt\\
 &=&\int_{-\pi}^{\pi} |f(e^{it})|^p|g(e^{it})|^pdt<\infty\}.
\eeqa
Clearly $H^p(|g|^p)=\{f\in \Hol(\D):fg\in H^p\}$, and
$f\lmto fg$ induces an isometry from $H^p(|g|^p)$ onto $H^p$.
Such spaces are not artificial. They appear naturally in the
context of Toeplitz operators. Indeed, if $\vp=\overline{Ig}/g$,
is a unimodular symbol, then $\ker T_{\vp}=gK^2_I(|g|^2)$
(see \cite{HS}). Here $T_{\vp}$ is defined in the usual
way by $T_{\vp}f=P_+(\vp f)$, $P_+$ being the standard
Riesz projection on $L^p(\T)$: $\sum_{n\in\Z} a_n\zeta^n
\lmto \sum_{n\ge 0}a_n\zeta^n$, $\zeta\in\T$.
Note that whenever $0\neq f\in \ker T_{\vp}$, where $\vp$ is
unimodular and $f=Jg$ is the inner-outer factorization of $f$,
then there exists an inner function $I$ such that $\vp=
\overline{Ig}/g$.

The representation $\ker T_{\vp}=gK^p_I(|g|^p)$ is particularly
interesting when $g$ is the extremal function of $\ker T_{\vp}$.
Then we know from a result by Hitt \cite{Hi} (see also
\cite{Sa94} for a de Branges-Rovnyak spaces approach to Hitt's result)
that when $p=2$, $\ker T_{\vp}=gK^2_I$, and that
$g$ is an isometric divisor on $\ker T_{\vp}=gK^2_I$
(or $g$ is an isometric multiplier on $K^2_I$).
In this situation we thus have $K^2_I(|g|^2)=K^2_I$.
Note, that for $p\neq 2$, 
if $g$ is extremal for
$gK^p_I(|g|^p)$, then $K^p_I(|g|^p)$ can
still be imbedded into $K^2_I$ when $p>2$ and in $K^p_I$ when
$p\in (1,2)$ (see \cite{HS}, where it is also shown that these
imbeddings can be strict).
In these situations when considering questions concerning
pseudocontinuation and analytic continuation, we can carry
over to $K^p_I(|g|^2)$
everything we know about $K^2_I$ or $K^p_I$ (which is the same
concerning these continuation matters).

However, in general the extremal function is not
easily detectable (explicit examples of extremal functions were
given in \cite{HS}), in that we cannot determine it, or for a
given $g$ it is not a simple matter to check whether it is extremal or not. 
So the first question that we would like to discuss is 
under which conditions on
$g$ and $I$, we can still say something about analytic continuation
of functions in $K^p_I(|g|^p)$. 
Our main result is that under a local integrability condition of $1/g$
on a closed arc not meeting the spectrum of $I$ it is possible to
extend every $K^p_I(|g|^p)$ function through such an arc. The
integrability condition is realized if for example $|g|^p$ is
an $(A_p)$ weight (but in this situation the analytic continuation 
turns out to be a simple consequence of H\"older's inequality and
Moeller's original result, see Proposition \ref{prop4n} and
comments thereafter).

In connection with analytic continuation under growth conditions,
another important result can be mentioned. Beurling 
(see \cite{Be}) proved that under some integral condition of
a weight $w$ defined on a square $Q$, {every} function holomorphic
on the upper and the lower half of the square and which is bounded
by a constant times $1/w$ extends analytically to the whole square.
Our result is different since we do not consider generic functions
holomorphic in both halfs of the square but admitting already
a certain type of pseudocontinuation. This allows to weaken the 
condition on the weight under which the analytic continuation 
is possible.

One could also
ask when $K^p_I(|g|^p)$ still embeds into $H^p$, or even in 
$H^r$, $r<p$, in other words,
when $K^p_I(|g|^p)\subset K^p_I$ or $K^p_I(|g|^p)\subset K^r_I$ 
(note that $K^p_I\cap \Hi
\subset K^p_I(|g|^p)$, which in particular gives
$K^p_I\subset K^p_I(|g|^p)$ whenever $g$ is bounded, so that
in such a situation the preceding inclusion is in fact an equality).
We will discuss some examples in this direction related to
our main result.

Naturally related to the question of analytic continuation
is the spectrum of the restriction of the backward shift
operator to $K^p_I(|g|^p)$ (see \cite{ARR}).
As was done in Moeller's paper,
we will explore these relations in the proof of our
main theorem on analytic continuations.

Finally we mention a paper by Dyakonov (\cite{dyak}).
He discussed Bernstein type inequalities in kernels
of Toeplitz operators which we know from our previous discussions
are closely related with weighted spaces $K^p_I(|g|^p)$.
More precisely, he discussed the regularity of functions in
$\ker T_{\vp}$ depending on the smoothness of the symbol $\vp$.\\

The paper is organized as follows. In the next section we will
discuss the analytic continuation when the spectrum of $I$ is
far from points where $g$ vanishes essentially. We will also 
establish a link with the spectrum of the backward shift on
$K^p_I(|g|^p)$ in this situation. As a corollary we deduce
that in certain situations one can get an embedding of
$K^p_I(|g|^p)$ into its unweighted companion. 
A simple situation is discussed when $K^p_I(|g|^p)$ can be 
embedded into a bigger $K^r_I$ ($1<r<p$),
and so still guaranteeing the analytic continuation
outside the spectrum of $I$.
Section 3 is different in flavour.
We will focus on the embedding problem by
discussing
some examples when $K^2_iI(|g|^2)$ does not embed into $K^2_I$.
It turns out that in the examples considered
the analytic continuation is like in the unweighted case.
Also, in these examples the spectrum of $I$ comes close to
points where $g$ vanishes essentially.

\section{Results when $\sigma(I)$ is far from the points where $g$ 
vanishes}\label{S2}

We start with a simple example.
Let $I$ be arbitrary with $-1\notin\sigma(I)$, 
and let $g(z)=1+z$, 
so that $\sigma(I)$ is
far from the only point where $g$ vanishes.
We know that $\ker T_{\frac{\overline{Ig}}{g}}=gK^p_I(|g|^p)$
(note that this is a so-called 
nearly invariant subspace). 
Let us compute this kernel. We first observe that 
$\frac{\overline{1+z}}{1+z}=\overline{z}$. 
Hence $T_{\frac{\overline{Ig}}{g}}
=T_{\overline{zI}}$, the kernel of which 
is known to be $K_{zI}^p$. So
\beqa
 gK^p_I(|g|^p)=K^p_{zI}.
\eeqa
The space on the right hand side contains the constant
functions. Hence $1/g\in K^p_I(|g|^p)$ 
(observe that 
$1/g=I\overline{\psi}$ with
$\psi(z)=Iz/(1+z)=Iz/g\in H^p_0(|g|^p)$).
In particular, $K^p_I(|g|^p)$ contains the function $1/g$ which is badly
behaved in $-1$, and thus cannot extend analytically through $-1$. 

This observation can be made more generally as stated in
the following result.

\begin{proposition}\label{prop1}
Let $g$ be an outer function in $H^p$. If
$\ker T_{\overline{g}/g}\neq\{0\}$ contains an {\rm inner}
function, then $1/g\in K^p_I(|g|^2)$ for \emph{every} inner
function $I$.
\end{proposition}

Before proving this simple result, we will do a certain number
of comments:

\begin{enumerate}
\item 
The example that we discussed above corresponds obviously to the situation
when the inner function contained in the kernel of 
$T_{\overline{g}/g}$ is identically equal to $1$.

\item The proposition again indicates a way of finding examples
of $g$ and $I$ such that $K^p_I(|g|^2)$ contains functions
that cannot be analytically continued through points where
$g$ is ``small''.

\item Suppose for the next two remarks that $p=2$.

\begin{itemize}
\item The claim that the kernel of $T_{\overline{g}/g}$ contains an
inner function implies in particular that this Toeplitz
operator is not injective and so $g^2$ is not rigid in $H^1$
(see \cite[X-2]{Sa}),
which means that
it is  not uniquely determined --- up to a real multiple --- 
by its argument (or equivalently, its normalized version $g^2/\|g^2\|_1$
is not exposed in the unit ball of $H^1$). 

\item It is clear that if the kernel of a Toeplitz operator is not
reduced to $\{0\}$ --- or equivalently (since $p=2$) $g^2$ is not rigid --- 
then it contains an {\it outer} function (just divide out the inner
factor of any non zero function contained in the kernel). 
However,
Toeplitz operators with non trivial kernels containing no inner
functions can be easily constructed. One could appeal to 
Hitt \cite{Hi}, Hayashi \cite{Hay} and Sarason \cite{Sa94}: 
the first states that every 
nearly invariant subspace is of the form $gK^2_I$ and $g$ is extremal
thus multiplying isometrically on $K^2_I$, the second tells us how
$g$ has to be in order that $gK^2_I$ is in particular the kernel of a
Toeplitz operator (the symbol being $\overline{Ig}/g$), and the
third one gives a general form of $g$ insuring the extremality
(or the isometric multiplication property). This  allows
to construct a kernel with the desired properties. 
However we can short-circuit these
results and take $T_{\overline{zg_0}/g_0}=T_{\overline{z}}
T_{\overline{g_0}/g_0}$, where $g_0(z)=(1-z)^{\alpha}$ and
$\alpha\in (0,1/2)$. The Toeplitz operator $T_{\overline{g_0}/g_0}$
is invertible ($|g_0|^2$ satisfies the Muckenhoupt $(A_2)$ condition --- see
Section \ref{S3} for more
discussions on invertibility of Toeplitz operators),
and 
$(T_{\overline{g_0}/g_0})^{-1}
=g_0P_+\frac{1}{\overline{g_0}}$ \cite{Ro} so that the kernel of
$T_{\overline{zg_0}/g_0}$ is given by the preimage under
$T_{\overline{g_0}/g_0}$ of the constants (which define the 
kernel of $T_{\overline{z}}$). Since $g_0P_+(c/\overline{g_0})
=cg_0/\overline{g_0(0)}$, $c$ being any complex number, we have
$\ker T_{\overline{zg_0}/g_0}=\C g_0$ which does not contain any inner
function.
\end{itemize}

\item 
See \cite{Ka} for a discussion of the intersection
$I_g:=H^p(|g|^p)\cap \overline{H_0^p(|g|^p)}$. 
Theorem 3 of that paper states that for points in the spectrum
of the inverse of the backward shift --- which is related with
the complement of those points where every 
function in $I_g$ continues analytically --- 
there always exist functions
with singularities in such a point.
\end{enumerate}

\begin{proof}[Proof of Proposition \ref{prop1}]
If $J\in \ker T_{\overline{g}/g}$, then
$\overline{g}/g=\overline{J\psi}$ where $\psi \in H^p_0$.
Since the functions appearing on both sides of the equality
are of modulus $1$, the function 
$\Theta=J\psi$ is inner and vanishes in 0. So 
\beqa
 T_{\overline{Ig}/g}=T_{\overline{I\Theta}}.
\eeqa
Hence 
\beqa
 K^p_I(|g|^p)=\frac{1}{g}\ker T_{\overline{I\Theta}}
 =\frac{1}{g}K^p_{I\Theta}.
\eeqa
Since $\Theta(0)=0$, we have $1\in K^p_{I\Theta}$, and so
$1/g\in K^p_I(|g|^2)$.
\end{proof}

One can also observe that if the inner function $J$
is in $\ker T_{\overline{g}/g}$
then $T_{\overline{Jg}/g}1=0$, and hence $1\in \ker
T_{\overline{Jg}/g}=gK^p_J(|g|^2)$ and $1/g\in K^p_J(|g|^2)$,
which shows that with this simpler argument the proposition 
holds for $I=J$.
Yet, our proof above allows to choose for $I$ any arbitrary inner
function and not necessarily that contained in $\ker T_{\overline{g}/g}$.
\\

So, without any condition on $g$, we cannot hope
for reasonable results. In the above example, when $p=2$,
then the function $g^2(z)=(1+z)^2$
is in fact not rigid (for instance the argument
of $(1+z)^2$ is the same as that of $z$).
Recall that
rigidity of $g^2$ is also characterized by the
fact that $T_{\overline{g}/g}$ is injective (see \cite[X-2]{Sa}). 
Here $T_{\overline{g}/g}=T_{\overline{z}}$ the kernel of which
is $\C$. From this it can also be deduced that $g^2$ is rigid
if and only if $H^p(|g|^p)\cap \overline{H^p(|g|^p)}=\{0\}$
which indicates again that rigidity should be assumed if we
want to have $K^p_I(|g|^p)$ reasonably defined.
\\

A stronger condition than rigidity (at least when $p=2$) is that of 
a Muckenhoupt weight.
Let us recall 
the Muckenhoupt $(A_p)$ condition: for general $1<p<\infty$ a
weight $w$ satisfies the $(A_p)$ condition if
\beqa
 B:=\sup_{I \text{ subarc of }\T}
 \left\{\frac{1}{|I|}\int_Iw(x)dx \times \left(\frac{1}{|I|}
 \int_Iw^{-1/(p-1)}(x)dx\right)^{p-1}\right\}<\infty.
\eeqa
When $p=2$, it is known that this condition is equivalent to
the so-called Helson-Szeg\H{o} condition.
The Muckenhoupt condition will play some r\^ole in the results
to come. However, our main theorem on analytic continuation (Theorem
\ref{thm1}) works under a weaker local integrability condition
which follows for instance from the Muckenhoupt condition.

Another observation can be made now.
We have already mentioned that the rigidity of $g^2$ in $H^1$
is equivalent to the injectivity of $T_{\overline{g}/g}$, when
$g$ is outer. It is also clear that $T_{g/\overline{g}}$ is
{\it always} injectif so that when $g^2$ is rigid, the operator
$T_{\overline{g}/g}$ is injectif with dense range. On
the other hand, by a result of Devinatz and Widom (see e.g.\
\cite[Theorem B4.3.1]{Nik}), the invertibility of $T_{\overline{g}/g}$,
where $g$ is outer, is equivalent to $|g|^2$ being $(A_2)$.
So the difference between rigidity and $(A_2)$ is the
surjectivity (in fact the closedness of the range)
of the corresponding Toeplitz operator. 
A criterion for surjectivity of non-injective Toeplitz
operators can be found in \cite{HSS}. 
It appeals to a parametrization which was earlier
used by Hayashi \cite{Hay} to characterize kernels of Toeplitz
operators among general nearly invariant subspaces. 
Rigid functions do appear in the characterization of Hayashi.

As a consequence of our main theorem (see Remark \ref{rem1})
analytic continuation can be expected on arcs not meeting the spectrum
of $I$ when $|g|^p$ is $(A_p)$. However the $(A_p)$
condition cannot be expected to be necessary since it is a global
condition whereas continuation depends on the local behaviour of
$I$ and $g$. We will
even give an example of a non-rigid function $g$
(hence not satisfying the $(A_p)$ condition) 
for which analytic continuation is always possible
in certain points of $\T$ where $g$ vanishes essentially. 
\\

Closely connected with backward shift invariant subspaces
is the spectrum of the
backward shift operator on the space under consideration.
The following result follows from 
\cite[Theorem 1.9]{ARR}: Let $B$ be the backward shift
on $H^p(|g|^p)$, defined by $Bf(z)=(f-f(0))/z$.
Clearly, $K^p_I(|g|^p)$ is invariant with respect to $B$
whenever $I$ is inner.
Then, 
$\sigma(B|K^p_I(|g|^p))
=\sigma_{ap}(B|K^p_I(|g|^p))$, where 
$\sigma_{ap}(T)=\{\lambda\in\C:\exists (f_n)_n$ with $\|f_n\|=1$
and $(\lambda-T)f_n\to 0\}$ denotes the
approximate point spectrum of $T$, and this
spectrum is equal to
\beqa
 \T\setminus\{1/\zeta\in\T:\mbox{ every }f\in K^p_I(|g|^p)
 \text{ extends analytically in a neighbourhood of }\zeta\}.
\eeqa

We would like to establish a link between this set and
$\sigma(I)$. To this end 
we will adapt the proof of the unweighted case 
\cite[Theorem 2.3]{Mo} to our situation.
As in the unweighted situation --- provided
the Muckenhoupt condition holds --- the approximate spectrum 
of $B|K^p_I(|g|^p)$ on $\T$ contains the conjugated spectrum 
of $I$. We will see later that the containment in the following
proposition actually
is an equality.

\begin{proposition}\label{prop3}
Let $g$ be outer in $H^p$ such that $|g|^p$ is a 
Muckenhoupt $(A_p)$-weight. Let $I$ be an inner function
with spectrum $\sigma(I)=\{\lambda\in \clos \D:\liminf_{z\to \lambda}I(z)=0\}$.
Then $\overline{\sigma(I)}\subset \sigma_{ap}(B|K^p_I(|g|^p))$.
\end{proposition}

\begin{proof}
It is clear that when $\lambda\in \D\cap\sigma(I)$, then 
$\lambda$ is a zero of $I$ and $k_{\lambda}^I=k_{\lambda}$.
It is a general fact that
$Bk_{\lambda}=\overline{\lambda}k_{\lambda}$. And since
$k_{\lambda}\in K^p_I\cap\Hi\subset K^p_I(|g|^p)$, we see
that $\overline{\lambda}$ is an eigenvalue, so it is in 
the point spectrum of $B|K^p_I(|g|^p)$ and then also in
the approximate point spectrum.

So, let us consider $\lambda\in\T$.
Take such a $\lambda\in \T$ with $\liminf_{z\to\overline{\lambda},z\in \D}
|I(z)|=0$. 
We want to prove that $\lambda\in
\sigma_{ap}(B|K^p_I(|g|^p))$.
To this end, let $\lambda_n\to \lambda$ a sequence of $\lambda_n\in\D$
with $I(\overline{\lambda_n})\to 0$.
Clearly
\beqa
 c_nk_{\overline{\lambda_n}}(z)= \frac{c_n}{1-\lambda_n z}
 =c_n\frac{1-\overline{I(\overline{\lambda_n})}I(z)}
 {1-\lambda_n z}+c_n
  \frac{\overline{I(\overline{\lambda_n})}I(z)}{1-\lambda_n z}
 =\underbrace{c_nk^I_{\overline{\lambda_n}}(z)}_{=:l_n(z)}
 +\underbrace{c_n\overline{I(\overline{\lambda_n})}
  I(z)k_{\overline{\lambda_n}}(z)}_{=:r_n(z)},
\eeqa
where $c_n$ is chosen so that $\|c_nk_{\overline{\lambda_n}}\|_{|g|^p}=1$.
Clearly $l_n=c_nP_Ik_{\overline{\lambda_n}}$ as a 
projected
reproducing kernel is in $K^p_I\cap\Hi$  which is contained 
in $K^p_I(|g|^p)$,
and $r_n\in I(H^p\cap \Hi)\subset IH^p(|g|^p)$. Since $|g|^p$ is
Muckenhoupt $(A_p)$, the Riesz projection $P_+$ is continuous on $L^p(|g|^p)$
and so also $P_I=IP_-\overline{I}$. As in 
the unweighted case $K^p_I(|g|^p)=P_IH^p(|g|^p)$ and $IH^p(|g|^p)=
\ker P_I|H^p(|g|^p)$, so that
the norm of ${c_n}k_{\overline{\lambda_n}}$ in $H^p(|g|^p)$
is comparable to the sum of the norms of $l_n$ and $r_n$:
\beqa
 1=\|{c_n}k_{\overline{\lambda_n}}\|_{|g|^p}
 \simeq \|l_n\|_{|g|^p}+\|r_n\|_{|g|^p}
 =\|l_n\|_{|g|^p}+|I(\overline{\lambda_n})|.
\eeqa
(For $p=2$, the boundedness of the projection $P_I$ means that
the angle between $K^2_I(|g|^2)$ and $IH^2(|g|^2)$ 
is strictly positive).
Since $I(\overline{\lambda_n})$ goes to zero, this implies that
the norms $\|l_n\|_{|g|^p}$ are comparable to a strictly positive 
constant. 
Now, 
\beqa
 \|(B-\lambda){c_n}k_{\overline{\lambda_n}}\|_{|g|^p}
 =\|c_n(\lambda_n-\lambda)k_{\overline{\lambda_n}}\|_{|g|^p}
 =|\lambda-\lambda_n|,
\eeqa
and hence
\beqa
 \|(B-\lambda)l_n\|_{|g|^p}
 &=&\|(B-\lambda){c_n}k_{\overline{\lambda_n}}
  - (B-\lambda)r_n\|_{|g|^p}\\
 &\le& |\lambda-\lambda_n|+\|B-\lambda\|_{H^p(|g|^p)\to H^p(|g|^p)}
  |I(\overline{\lambda_n})|\\
 &\le& |\lambda-\lambda_n|+(\|B\|_{H^p(|g|^p)\to H^p(|g|^p)}+|\lambda|)
 |I(\overline{\lambda_n})|,
\eeqa
which tends to zero, while $\|l_n\|_{|g|^p}$ is uniformly bounded
away from zero. So, $\lambda\in \sigma_{ap}(B|K^p_I(|g|^p)$.
\end{proof}

We now come to our main theorem.

\begin{theorem}\label{thm1}
Let $g$ be an outer function in $H^p$, $1<p<\infty$ and
$I$  an inner function with associated
spectrum $\sigma(I)$. Let $\Gamma$ be a
closed arc in $\T$ not meeting
$\sigma(I)$. If there exists
$s>q$, $\frac{1}{p}+\frac{1}{q}=1$, with $1/g\in L^s(\Gamma)$, then
every function $f\in K^p_I(|g|^p)$ extends
analytically through $\Gamma$.
\end{theorem} 

\begin{remark}\label{rem1}
It is known (see e.g. \cite{Mu}) that when $|g|^p\in (A_p)$, $1<p<\infty$,
then there exists $r_0\in (1,p)$ such that $|g|^p\in (A_r)$ for
every $r>r_0$. Take $r\in (r_0,p)$.
Then in particular $1/g\in L^{s}$, where
$\frac{1}{r}+\frac{1}{s}=1$. Since $r<p$ we have $s>q$.
which allows to conclude that in this situation
$1/g\in L^s(\Gamma)$ for every $\Gamma\subset\T$ ($s$ independant
of $\Gamma$).
\end{remark}

We promised earlier an example of a non-rigid function $g$
for which analytic continuation of $K^p_I$-functions
{\it is} possible in certain points where $g$ vanishes. 

\begin{example}
For $\alpha\in (0,1/2)$, let  $g(z)=(1+z)(1-z)^{\alpha}$.
Clearly $g$ is an outer function vanishing essentially in $1$ and $-1$. 
Set $h(z)=z(1-z)^{2\alpha}$, then
by similar arguments as those employed in the introducing example to this
section one can check that $\arg g^2=\arg h$ a.e.\ on $\T$. Hence
$g$ is not rigid (it is the ``big'' zero in $-1$ which is
responsible for non-rigidity). On the other hand, the zero in $+1$
is ``small'' in the sense that $g$ satisfies the local integrability
condition in a neighbourhood of $1$
as required in the theorem, so that whenever $I$ has its
spectrum far from $1$, then every $K^2_I(|g|^2)$-function can be analytically
continued through suitable arcs around $1$.

This second example can be pushed a little bit further. In the spirit of
Proposition \ref{prop1} we check that (even) when
the spectrum of an inner function $I$ does not meet $-1$, there are
functions in $K^p_I(|g|^p)$ that are badly behaved in $-1$.
Let again $g_0(z)=(1-z)^{\alpha}$. Then
\beqa
 \frac{\overline{g(z)}}{g(z)}=\frac{\overline{(1+z)(1-z)^{\alpha}}}
 {(1+z)(1-z)^{\alpha}}=\overline{z}\frac{\overline{g_0(z)}}{g_0(z)}.
\eeqa
As already explained, for every inner function $I$, we have
$\ker T_{\overline{Ig}/g}=gK^p_I(|g|^p)$, so that we are
interested in the kernel $\ker T_{\overline{Ig}/g}$.
We have $T_{\overline{Ig}/g}f=0$ when $f=Iu$ and
$u\in \ker T_{\overline{g}/g}=\ker T_{\overline{zg_0}/g_0}
=\C g_0$ (see the discussion just
before the proof of Proposition \ref{prop1}).
Hence the function defined by
\beqa
 F(z)=\frac{f(z)}{g(z)}={I(z)} 
 \frac{g_0(z)}{g(z)}=
 \frac{I(z)}{1+z}
\eeqa
is in $K^p_I(|g|^p)$ and it is badly behaved in $-1$ when the
spectrum of $I$ 
does not meet $-1$ (but not only). 
\end{example}

The preceding discussions motivate the following question: does
rigidity of $g$ suffice to get analytic continuation for 
$K^p_I(|g|^2)$-function whenever $\sigma(I)$ is far from 
zeros of $g$?

\begin{proof}[Proof of Theorem \ref{thm1}]
Take an arc as in the theorem. 
Since $\sigma(I)$ is closed, the distance between
$\sigma(I)$ and $\Gamma$ is strictly positif.
Then there is a neigbourhood of $\Gamma$ intersected with $\D$
where $|I(z)|\ge \delta>0$. It is clear that in this neighbourhood we
are far away from the spectrum of $I$.
Thus  $I$ extends 
analytically through $\Gamma$.
For
what follows we will call the endpoints of this arc $\zeta_1$
and $\zeta_2$.
We would like to know whether every function in $K^p_I(|g|^p)$
extends also analytically through $\Gamma$. 

We adapt an argument by Moeller based on Morera's theorem.
Let us first introduce some notation (see Figure 1). 

\begin{center}
\begin{picture}(0,0)%
\includegraphics{sector.pstex}%
\end{picture}%
\setlength{\unitlength}{1579sp}%
\begingroup\makeatletter\ifx\SetFigFont\undefined%
\gdef\SetFigFont#1#2#3#4#5{%
  \reset@font\fontsize{#1}{#2pt}%
  \fontfamily{#3}\fontseries{#4}\fontshape{#5}%
  \selectfont}%
\fi\endgroup%
\begin{picture}(8274,6924)(1789,-7573)
\put(7051,-3286){\makebox(0,0)[lb]{\smash{{\SetFigFont{6}{7.2}{\familydefault}{\mddefault}{\updefault}{\color[rgb]{0,0,0}$\zeta_1$}%
}}}}
\put(6301,-2386){\makebox(0,0)[lb]{\smash{{\SetFigFont{6}{7.2}{\familydefault}{\mddefault}{\updefault}{\color[rgb]{0,0,0}$\zeta_2$}%
}}}}
\put(6976,-5761){\makebox(0,0)[lb]{\smash{{\SetFigFont{8}{9.6}{\familydefault}{\mddefault}{\updefault}{\color[rgb]{0,0,0}$\T$}%
}}}}
\put(6151,-4486){\makebox(0,0)[lb]{\smash{{\SetFigFont{6}{7.2}{\familydefault}{\mddefault}{\updefault}{\color[rgb]{0,0,0}$r_1$}%
}}}}
\put(6751,-4486){\makebox(0,0)[lb]{\smash{{\SetFigFont{6}{7.2}{\familydefault}{\mddefault}{\updefault}{\color[rgb]{0,0,0}$r_2$}%
}}}}
\put(7426,-4486){\makebox(0,0)[lb]{\smash{{\SetFigFont{6}{7.2}{\familydefault}{\mddefault}{\updefault}{\color[rgb]{0,0,0}$1/r_2$}%
}}}}
\put(8551,-4486){\makebox(0,0)[lb]{\smash{{\SetFigFont{6}{7.2}{\familydefault}{\mddefault}{\updefault}{\color[rgb]{0,0,0}$1/r_1$}%
}}}}
\put(5101,-2611){\makebox(0,0)[lb]{\smash{{\SetFigFont{6}{7.2}{\familydefault}{\mddefault}{\updefault}{\color[rgb]{0,0,0}$C(r_2)$}%
}}}}
\put(6451,-1561){\makebox(0,0)[lb]{\smash{{\SetFigFont{6}{7.2}{\familydefault}{\mddefault}{\updefault}{\color[rgb]{0,0,0}$\tilde{C}(r_2)$}%
}}}}
\end{picture}%
\\
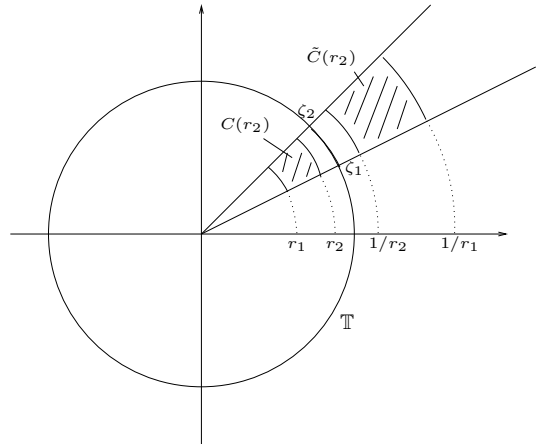
Figure 1: The regions $C(r_2)$ and $\tilde{C}(r_2)$ 
\end{center}

For suitable $r_0\in (0,1)$ let
$r_0<r_1<r_2<1$. 
Using Moeller's notation, we call $C(r_2)$ the boundary of a circular
sector whose vertices are $\{r_k\zeta_l\}_{k,l=1,2}$ (inside $\D$)
and $\tilde{C}(r_2)$ is the boundary of a circular sector
with vertices $\{r_k^{-1}\zeta_l\}_{k,l=1,2}$ (outside $\D$).
These sectors are thus in a sense symmetric with respect
to $\T$ and when $r_2$ goes to 1, then at the limit they will form
a circular sector with vertices $\{r_1\zeta_l\}_{l=1,2}$ and
$\{r_1^{-1}\zeta_l\}_{l=1,2}$ the interior of which contains 
in particular $\Gamma$. 
For the construction to come, we need to assume that $r_1$ 
is chosen in such a way that $1/I$ is analytic (and bounded)
in $C(r_2)$ for every $r_1<r_2<1$, which is of course possible.
We will also assume that $f$ admits boundary values in $\zeta_1$
and $\zeta_2$. Since $f$ is in $N^+$ it has boundary values
a.e., and so the previous requirement is not difficult to meet
(by possibly moving the endpoints $\zeta_1,\zeta_2$ if necessary).

Take now $f\in K^p_I(|g|^p)=H^p(|g|^p)
\cap I\overline{H^p_0(|g|^p)}$, so that $f=I\overline{\psi}$,
where $\psi\in H^p_0(|g|^p)$ can be written as
$\psi(z)=\sum_{n\ge 1} c_nz^n$, $z\in\D$. 
The function $\tilde{\psi}$ defined by
\beqa
 \tilde{\psi}(z)=\sum_{n\ge 1}\overline{c_n}\frac{1}{z^n},\quad
 z\in \D_e:=\hat{C}\setminus \clos \D,
\eeqa
(the tilde-sign does not mean harmonic conjugation here)
is the pseudocontinuation of the meromorphic function
$f/I$. Note that $f/I$ is even analytic in 
every $C(r_2)$, $r_1<r_2<1$. 
Since $f\in K^p_I(|g|^p)$, we in particular have
\beqa
 \sup_{r<1}\int_{\Gamma}|f(re^{it})|^p|g(re^{it})|^p\frac{dt}{2\pi}<\infty.
\eeqa
Then, since by construction $1/I$ is analytic and bounded on 
$C(r_2)$, $r_0<r_1<r_2<1$, we also have
\beqa
 M_1:=\sup_{r_1<r<1}\int_{\Gamma}|
 \frac{f(re^{it})}{I(re^{it})}|^p|g(re^{it})|^p
 \frac{dt}{2\pi}<\infty.
\eeqa
On the other hand,
for $r_1<r<1$, we have
\beqa
 \int_{\Gamma}|\tilde{\psi}\left(\frac{e^{it}}{r}\right)|^p
 |g(re^{it})|^p\frac{dt}{2\pi}
 &=&\int_{\Gamma}|\psi(re^{it})|^p
 |g(re^{it})|^p\frac{dt}{2\pi}\\
 &\le& \sup_{r_1<\rho<1} \int_{\Gamma}|\psi(\rho e^{it})|^p
 |g(\rho e^{it})|^p\frac{dt}{2\pi}=:M_2<\infty
\eeqa
since $\psi\in H^p_0(|g|^p)$.
We will now use a version of Lebesgue's monotone convergence theorem.
Recall that 
$\tilde{\psi}$ is the
pseudocontinuation of $f/I$, so that
a.e.\ on $\T$ and in particular on $\Gamma$,
\bea\label{eq6}
 \frac{f(re^{it})}{I(re^{it})}-\tilde{\psi}(\frac{e^{it}}{r})\to 0,
 \quad \text{ when }r\to 1^-.
\eea

Let us prove that there exists $\eps>0$ and $M>0$ such that
for every $r_1<r<1$,
\bea\label{eq4}
 \int_{\Gamma} \left|\frac{f(re^{it})}{I(re^{it})}-\tilde{\psi}
 (\frac{e^{it}}{r})\right|^{1+\eps}\frac{dt}{2\pi}\le M.
\eea
Note first that for $r_1<r<1$
\beqa
 \lefteqn{ \int_{\Gamma} \left|\frac{f(re^{it})}{I(re^{it})}-\tilde{\psi}
 (\frac{e^{it}}{r}) \right|^{p}|g(re^{it})|^p\frac{dt}{2\pi}}\\
 &&\le c_p\left(\int_{\Gamma} \left|\frac{f(re^{it})}{I(re^{it})}\right|^{p}
    |g(re^{it})|^p\frac{dt}{2\pi}
 +\int_{\Gamma} \left|\tilde{\psi}
 (\frac{e^{it}}{r}) \right|^{p}|g(re^{it})|^p\frac{dt}{2\pi}\right)\\
 &&\le c_p(M_1+M_2)=:M'. 
\eeqa
From this and the 
condition $1/g\in L^s(\Gamma)$ we will 
deduce via a simple H\"older inequality our estimate \eqref{eq4}.

By assumption there is $s>q$ with $1/g\in L^s(\Gamma)$. Then there exists
$\eps>0$ such that $s>p\frac{1+\eps}{p-(1+\eps)}>q=p/(p-1)$. Hence
$L^s(\Gamma)\subset L^{p\frac{1+\eps}{p-(1+\eps)}}(\Gamma)$ and
\bea\label{eq13}
 \int_{\Gamma}\left(\frac{1}{|g|}\right)^
 {p\frac{1+\eps}{p-(1+\eps)}}dm
 \le c \int_{\Gamma}\frac{1}{|g|^s}dm<\infty
\eea

Then by H\"older's inequality, assuming also $p/(1+\eps)>1$,
\beqa
 \lefteqn{\int_{\Gamma} \left|\frac{f(re^{it})}{I(re^{it})}-\tilde{\psi}
 (\frac{e^{it}}{r})\right|^{1+\eps}\frac{dt}{2\pi}
 =\int_{\Gamma} \left|\frac{f(re^{it})}{I(re^{it})}-\tilde{\psi}
 (\frac{e^{it}}{r})\right|^{1+\eps}\frac{|g(re^{it})|^{1+\eps}}
 {|g(re^{it})|^{1+\eps}}\frac{dt}{2\pi}}\\
 &&\le \left\{ \int_{\Gamma} \left|\frac{f(re^{it})}{I(re^{it})}-\tilde{\psi}
 (\frac{e^{it}}{r})\right|^p|g(re^{it})|^p \frac{dt}{2\pi}\right\}^{(1+\eps)/p}
   \left\{ \int_{\Gamma} \left(\frac{1}{|g(re^{it})|^{1+\eps}}
  \right)^{p/(p-(1+\eps))}\frac{dt}{2\pi}\right\}^{(p-(1+\eps))/p}
\eeqa
The first factor in this product is uniformly bounded by our
previous discussions. Consider the second factor. Recall that $g$
is outer (hence $1/g$ is in the Smirnov class) so that it is 
sufficient to check whether the second factor is bounded for
$r=1$, and this follows from \eqref{eq13}.
We have thus proved \eqref{eq4}.

By standard arguments based on Lebesgue's dominated convergence
theorem and Tchebychev's inequality we get  
that \eqref{eq4} together with \eqref{eq6} imply that
\beqa
 \int_{\Gamma} \frac{f(re^{it})}{I(re^{it})}
   -\tilde{\psi}(\frac{e^{it}}{r})\frac{dt}{2\pi}\to 0.
\eeqa

From this point on we can repeat Moeller's proof \cite[Lemma 2.2]{Mo},
which is based on Morera's theorem, and which uses the fact that
when $r_2\to 1^-$, then the regions $C(r_2)$ and $\tilde{C}(r_2)$
fusion to a big angular sector, where the arcs $\Gamma_{r_2}
:=\{r_2\zeta:\zeta\in\Gamma\}$ and
$\Gamma_{1/r_2}$ are oriented in an opposite direction and the
difference of the integrals of our fonctions $f/I$ 
and $\psi$ on these two arcs tends to zero. 

This implies that $f/I$ extends analytically through $\Gamma$.
\end{proof}

The theorem together with Proposition \ref{prop3} and Remark
\ref{rem1} allow us to obtain the following result.

\begin{corollary}\label{cor1}
Let $g$ be outer in $H^p$ such that $|g|^p$ is a 
Muckenhoupt $(A_p)$ weight. Let $I$ be an inner function
with spectrum $\sigma(I)=\{\lambda\in \clos\D:\liminf_{z\to \lambda}I(z)=0\}$.
Then $\overline{\sigma(I)}=\sigma_{ap}(B|K^p_I(|g|^p)$.
\end{corollary}

\begin{proof}
The inclusion $\overline{\sigma(I)}\subset\sigma_{ap}(B|K^p_I(|g|^p))$
has been discussed in Proposition \ref{prop3}.

For the reverse inclusion, suppose that $\overline{\lambda}\not\in
\sigma(I)$. 
Note that $\overline{\sigma(I)}\cap \D=\sigma_p(B|K^p_I(|g|^p))
=\sigma_{ap}(B|K^p_I(|g|^p))=\sigma(B|K^p_I(|g|^p))$ so that it
is sufficient to consider the case $\lambda\in\T$.
Since $\sigma(I)$ is closed, there is an arc not meeting
$\sigma(I)$ and containing $\overline{\lambda}$. By Theorem
\ref{thm1} and Remark \ref{rem1}, 
every $f\in K^p_I(|g|^p)$ extends analytically through
this arc, and thus in particular through $\overline{\lambda}$.
So, by \cite[Theorem 1.9]{ARR}, $\overline{\lambda}$ cannot
be in $\sigma_{ap}(B|K^p_I(|g|^p))$ (neither in $\sigma(B|K^p_I(|g|^p))$).
\end{proof}

Another simple consequence of Theorem \ref{thm1} concerns embeddings.

\begin{corollary}\label{cor2}
Let $I$ be an inner function with spectrum $\sigma(I)$.
If $\Gamma\subset \T$ is a closed arc not meeting $\sigma(I)$ 
and if $g$ is an outer function in $H^p$ such that
$|g|\ge \delta$ on $\T\setminus \Gamma$ for some constant $\delta>0$
and $1/g\in L^s(\Gamma)$, $s>q$, $\frac{1}{p}+\frac{1}{q}=1$. 
Then $K^p_I(|g|^p)\subset K^p_I$.
If moreover $g$ is bounded, then the last inclusion is an
equality.
\end{corollary}

\begin{remark}
1) Suppose $p=2$.
By Hitt's result \cite{Hi}, when $g$ is
the extremal function of a nearly invariant subspace $M\subset H^2$, 
then there exists an inner function 
$I$ such that $M=gK^2_I$, and $g$ is an isometric multiplier
on $K^2_I$ so that $K^2_I=K^2_I(|g|^2)$. With the corollary we
can construct examples of spaces of infinite dimension
where the latter identity holds without
$g$ being extremal. Recall from \cite[Lemma 3]{HS} that a
function $g$ is extremal for $gK^2_I(|g|^2)$ if
$\int f |g|^2dm=f(0)$ for every function $f\in K^2_I(|g|^2)$. 
The following example is constructed in the
spirit of the example in \cite[p.356]{HS}. 
Fix $\alpha\in (0,1/2)$. Let $\gamma(z)=(1-z)^{\alpha}$ 
and let $g$ be an outer
function in $H^2$ such that $|g|^2=\Re \gamma$
a.e.\ on $\T$
(such a function clearly exists). 
Let now $I=B_{\Lambda}$ be a Blaschke product with $0\in\Lambda$.
If $\Lambda$ accumulates to points outside $1$, then the corollary
shows that $K^2_I=K^2_I(|g|^2)$. Let us check that $g$ is
not extremal. To this end we compute $\int k_{\lambda}|g|^2dm$
for $\lambda\in\Lambda$ (recall that for $\lambda\in\Lambda$,
$k_{\lambda}\in K^2_I=K^2_I(|g|^2)$):
\beqa
 \int k_{\lambda}|g|^2dm
 &=&\int k_{\lambda}\Re \gamma dm
 =\frac{1}{2}\left(\int k_{\lambda}\gamma dm
  +\int k_{\lambda}\overline{\gamma}dm\right)  
 =\frac{1}{2}k_{\lambda}(0)\gamma(0)+\frac{1}{2}\langle
 k_{\lambda},\gamma\rangle\\
 &=&\frac{1}{2}(1+\overline{(1-\lambda)^{\alpha}})
\eeqa
which is different from $k_{\lambda}(0)=1$ (except when $\lambda=0$). 
Hence $g$ is not extremal.

2) Observe that {\bf if} $K^p_I(|g|^p)\subset K^p_I$, then the
analytic continuation is of course a simple consequence of that 
in $K^p_I$ (and hence of Moeller's result). And since 
$K^p_I(|g|^2)$ always contains $k^I_{\lambda}$ which continues
only outside $\sigma(I)$, one cannot hope
for a better result. Note also that the 
inclusion $K^p_I(|g|^p)\subset K^p_I$ is
a kind of reverse inclusion to those occuring in the context of
Carleson measures. Indeed, the problem of knowing when
$K^p_I\subset L^p(\mu)$ continues attracting a lot of attention 
(and the reverse situation to ours would correspond to
$d\mu=|g|^2dm$). Such a measure $\mu$ is called a Carleson measure
for $K^p_I$, and it is notoriously difficult to describe these
in the general situation (see \cite{TV1} for some results;
when $I$ is a so-called one-component inner function $I$, a
geometric characterization is available).
\end{remark}

\begin{proof}[Proof of Corollary \ref{cor2}]
Pick $f\in K^p_I(|g|^2)$. We only have to prove that $f\in L^p$.
By Theorem \ref{thm1}, $f$ continues analytically through
$I$ and so $f$ is bounded on $I$. On the other hand
\beqa
 \int_{\T\setminus \Gamma}|f|^pdm\le \frac{1}{\delta^p}
 \int_{\T\setminus \Gamma} |f|^p|g|^pdm<\infty.
\eeqa
\end{proof}

It is clear that the corollary is still valid when $\Gamma$
is replaced by a finite union of intervals.
However, in the next section we will see that it is no longer valid
when $\Gamma$ is replaced by an infinite union of intervals
under a yet weaker integrability condition on $1/g$ (see
Proposition \ref{prop5n}).

A final simple observation concerning the local integrability
condition $1/g\in L^s(\Gamma)$, $s>q$: if it is replaced by
the global condition $1/g\in L^s(\T)$, then we have an
embedding into a bigger backward shift invariant subspace:

\begin{proposition}\label{prop4n}
Let $1<p<\infty$ and $1/p+1/q=1$. If there exists $s>q$ such that
$1/g\in L^s(\T)$, then for $r$ with $1/r=1/p+1/s$
we have $L^p(|g|^p)\subset L^r$.
\end{proposition}

\begin{proof}
This is a simple application of Young's (or H\"older's) inequality:
\beqa
 \int |f|^rdm=\int |fg|^r\frac{1}{|g|^r}dm
 \le \left(\int |fg|^pdm\right)^{r/p}\left(\int\frac{1}{|g|^s}
 dm\right)^{r/s}<\infty.
\eeqa
\end{proof}

Under the assumptions of the proposition we of course also
have $K^p_I(|g|^p)\subset K^r_I$.
In particular, under the
assumption of the proposition, Moeller's theorem
then shows that every function $f\in K^p_I(|g|^p)$
extends analytically outside $\sigma(I)$.

Another observation is that when $|g|^p\in (A_p)$, then
by Remark \ref{rem1}, 
we have $1/g\in L^s$ for some $s>q$ and so 
the assumptions of the proposition are met, and
again $K^p_I(|g|^p)$ embeds into 
$K^{r}_I$.

\section{Examples when $K^2_I(|g|^2)$ does not embed into
$K^2_I$}
\label{S3}

Here we will discuss some examples when $K^2_I(|g|^2)$ does not
embed into $K^2_I$ even when $|g|^2$ satisfies some regularity
condition. The first example is when $|g|^2$ is $(A_2)$. The second
example, discussed in Proposition \ref{prop5n}, is a kind of
counterpart to Corollary \ref{cor2}. In both examples the spectrum
of $I$ comes close to the points where $g$ vanishes essentially. 

Before entering into the discussion of our examples, we
give a result in connection with
invertibility of Toeplitz operators.  
For an outer function $g\in H^p$,
the Toeplitz operator $T_{\overline{g}/g}$ is invertible
if and only if $|g|^p$ is an $(A_p)$ weight (we have
already mentioned the result by Devinatz and Widom for the case $p=2$,
see e.g. \cite[Theorem B4.3.1]{Nik}; for general $p$, see \cite{Ro}).
If this is the case,
the inverse of $T_{\overline{g}/g}$ is defined 
by $A=gP_+\frac{1}{\overline{g}}$ (see \cite{Ro}). Then, 
the operator $A_0=P_+\frac{1}{\overline{g}}$
is an isomorphism of $H^p$ onto $H^p(|g|^p)$. 

\begin{lemma}\label{L1}
Suppose $|g|^p$ is an $(A_p)$ weight and $I$ 
an inner function.
Then $A_0=P_+\frac{1}{\overline{g}}$
is an isomorphism of $K^p_I$ onto $K^p_I(|g|^p)$.
Also, for every $\lambda\in \D$ 
we have
\bea\label{eq9}
 A_0k_{\lambda}=\frac{k_{\lambda}(\mu)}{\overline{g(\lambda)}}.
\eea
\end{lemma}

\begin{proof}
Let us first discuss the action of $A_0$ on the reproducing kernels:
\beqa
 A_0k_{\lambda}(\mu)=(P_+\frac{k_{\lambda}}{\overline{g}})(\mu)
 =\langle \frac{k_{\lambda}}{\overline{g}},k_{\mu}\rangle
 =\langle k_{\lambda},\frac{k_{\mu}}{{g}}\rangle
 =\overline{\left(\frac{k_{\mu}(\lambda)}{g(\lambda)}\right)}
 =\frac{k_{\lambda}(\mu)}{\overline{g(\lambda)}},
\eeqa
so that $A_0k_{\lambda}=k_{\lambda}/\overline{g(\lambda)}$. 

Note that from this we can deduce the inclusion 
$A_0K^p_I\subset K^p_I(|g|^p)$
when $I$ is a Blaschke product with simple zeros $\Lambda$,
since in that case $k^I_{\lambda}=k_{\lambda}$, $\lambda\in\Lambda$,
span the space $K^p_I$, and $k_{\lambda}=k^I_{\lambda}\in H^p(|g|^p)\cap
I\overline{H^p_0(|g|^p)}=K^p_I(|g|^p)$.

For general inner functions $I$ we need a different argument.
Recall that $|g|^p$ is an $(A_p)$ weight. So, the function $G:=1/g$ is
in $H^q$, $1/p+1/q=1$. 
Taking Fej\'er polynomials $G_N$ of $G$, we get 
$G_N\to G$ in $H^q$ (e.g.\ \cite[A3.3.4]{Nik}), $G_N\in \Hi$. Then
$A_N:=P_+\overline{G}_N$ is a finite linear combination of
composed backward shifts thus leaving $K^q_I$ invariant, so that 
$P_+\overline{G}_Nk_{\lambda}^I\in K^q_I$
for every $\lambda\in \D$. 
Now, since
$\overline{G}_Nk_{\lambda}^I\to \overline{G}k_{\lambda}^I$
in $H^q$, we obtain $A_Nk_{\lambda}^I\to A_0k_{\lambda}^I$
which is thus in $K^q_I$. On the other hand, since $k_{\lambda}^I
\in \Hi\subset H^p$, we also have $A_0k_{\lambda}^I\in H^p(|g|^p)$.
Hence $A_0k_{\lambda}^I\in K_I^q\cap H^p(|g|^p)\subset K^p_I(|g|^p)$.
Note that the projected reproducing
kernels $k_{\lambda}^I$, $\lambda\in\D$, generate a dense subspace
in $K^p_I$, so $A_0K^p_I\subset K^p_I(|g|^p)$.

Let us prove that $A_0$ is from $K^p_I$ onto $K^p_I(|g|^p)$.
To this end, let $h\in K^p_I(|g|^p)$, then
$gh\in gK^p_I(|g|^p)=\ker T_{\overline{Ig}/g}$, and since
$T_{\overline{Ig}/g}=T_{\overline{I}}T_{\overline{g}/g}$ with
$T_{\overline{g}/g}$ invertible, we have $gh\in \ker
T_{\overline{Ig}/g}$ if and only if $T_{\overline{g}/g}(gh)\in \ker 
T_{\overline{I}}=K^2_I$. And so $gh\in T_{\overline{g}/g}^{-1}K^2_I
=gP_+\frac{1}{\overline{g}}K^p_I$, from where we get $h\in 
P_+\frac{1}{\overline{g}}K^p_I=A_0K^p_I$.
\end{proof}

We refer the reader to \cite{dyak2}, in particular Proposition
1.3, for some discussions
of the action of $T_{\overline{G}}$ on $K^p_I$ spaces.

\begin{proposition}
There exists an outer function $g$ in $H^2$ with 
$|g|^2$ Muckenhoupt $(A_2)$, and an inner function $I$ such
that $K^2_I(|g|^2)\not\subset K^2_I$.
\end{proposition}

\begin{proof}
Take $g(z)=(1-z)^{\alpha}$ with $\alpha\in (0,1/2)$. Then
$|g|^2$ is $(A_2)$. 
Let also $I=B_{\Lambda}$ where
$\Lambda=\{1-1/2^n\}_n$. In this situation, $\sigma(I)\cap \T
=\{1\}$, which corresponds to the point where $g$ vanishes.
Clearly, since $\Lambda$ is an interpolating sequence, the
sequence $\{k_{\lambda_n}/\|k_{\lambda_n}\|_2\}_n$ is a
normalized unconditional basis
in $K^2_I$. This means that we can write
$K^2_I=l^2(\frac{k_{\lambda_n}}{\|k_{\lambda_n}\|_2})$ meaning
that $f\in K^2_I$ if and only if
\beqa
 f= \sum_{n\ge 1}\alpha_n \frac{k_{\lambda_n}}{\|k_{\lambda_n}\|_2}
\eeqa
with $\sum_{n\ge 1} |\alpha_n|^2<\infty$ (the last sum defines
the square of an equivalent
norm in $K^2_I$).

So, since $|g|^2$ is 
Muckenhoupt $(A_2)$, we get from \eqref{eq9}
\beqa
 \{A_0(k_{\lambda_n}/\|k_{\lambda_n}\|_2)\}_n
 =\{\frac{k_{\lambda_n}}{\overline{g(\lambda_n)}
 \|k_{\lambda_n}\|_2}\}_n,
\eeqa
and $\{k_{\lambda_n}/(\overline{g(\lambda_n)}{\|k_{\lambda_n}\|_2})\}_n$
is thus an unconditional basis in $K^2_I(|g|^2)$
(almost normalized in the sense that
$\|A_0(k_{\lambda_n}/\|k_{\lambda_n}\|_2)\|_{|g|^2}$ 
is comparable
to a constant independant of $n$). Hence for every sequence 
$\alpha=(\alpha_n)_n$ with $\sum_{n\ge 1}|\alpha_n^2|<\infty$,
we have
\beqa
 f_{\alpha}:=\sum_{n\ge 1}\alpha_n\frac{k_{\lambda_n}}{\overline{g(\lambda_n)}
 \|k_{\lambda_n}\|_2}\in K^2_I(|g|^2).
\eeqa
Now, in order to construct a function in $K^2_I(|g|^2)$ which
is not in $K^2_I$, it suffices to choose $(\alpha_n)_n$ in
such a way that 
\beqa
 \sum_{n\ge 1}|\alpha_n|^2<+\infty,
\eeqa
so that $f_{\alpha}\in K^2_I(|g|^2)$, and
\beqa
 \sum_{n\ge 1}\left|\frac{\alpha_n}{\overline{g(\lambda_n)}}\right|^2
 =+\infty,
\eeqa
so that $f_{\alpha}\not\in K^2_I$.
Recall that $\lambda_n=1-1/2^n$ and $g(z)=(1-z)^{\alpha}$. Hence
$|g(\lambda_n)|^2=2^{-2n\alpha}$.
Now, taking e.g. $\alpha_n=1/n$ the first of the above two sums
converges, and $|\alpha_n/\overline{g(\lambda_n)}|^2=2^{2n\alpha}/n^2$
which does even not converge to zero.
\end{proof}

According to Proposition \ref{prop4n}, the function $f_{\alpha}$
constructed in the proof is in some $K^r_I$ for a suitable $r\in
(1,2)$ (this can also be checked directly by choosing $r\in (1,2)$ 
in such a way that $(\frac{\alpha_n}{\overline{g(\lambda_n)}}
\frac{\|k_{\lambda_n}\|_r}{\|k_{\lambda_n}\|_2})_n\in l^r$, which
is possible).\\

In view of Corollary \ref{cor2}, we will discuss another
situation. In that corollary we obtained that when $g$ is uniformly bounded
away from zero on $\T\setminus\Gamma$ where $\Gamma$ is an arc on
which $1/g$ is $s$-integrable, $s>q$, and $\Gamma$ not meeting
$\sigma(I)$, then $K^p_I(|g|^p)$ embeds into $K^p_I$ . 
We will now be interested in the situation
when the arc $\Gamma$ is replaced by an infinite union of arcs.
Our example is constructed for $p=2$.
Then under the weaker assumption $1/g\in L^s$, $s<2$, the
embedding turns out to be false in general.

\begin{proposition}\label{prop5n}
There exists an inner function $I$, a sequence of disjoint 
closed arcs $(\Gamma_n)_n$ in $\T$ not meeting the spectrum
of $I$, such that for every $s<2$ there is an outer function
$g\in H^2$ with $1/g\in L^s(\Gamma)$ 
and $|1/g|\ge \delta$ on $\T\setminus\Gamma$, where $\Gamma
=\bigcup_n\Gamma_n$, but $K^2_I(|g|^2)\not\subset K^2_I$.
\end{proposition}

\begin{proof}
Let $I(z)=\exp \frac{z+1}{z-1}$ which is a singular inner
function the associated measure of which is supported on $\{1\}$
(which is equal to $\sigma(I)$).
As in the preceding proposition we choose $\Lambda=\{\lambda_n\}_n
=\{1-1/2^n\}_n$, which is an interpolating sequence in 
$H^2$. Moreover, $I(\lambda_n)\to 0$ when $n\to\infty$, so that
$\Lambda$ is also an interpolating sequence for $K^2_I$
(see \cite{HNP} or \cite[D4.4.9 (8)]{Nik}). Set $f(z)=(1-I(z))/(1-z)$,
then 
\beqa
 |f(\lambda_n)|^2=\left|\frac{1-I(\lambda_n)}{1-\lambda_n}\right|^2
 =\frac{1-e^{1-2^{n+1}}}{1/2^{2n}}\simeq 2^{2n}
\eeqa
So, $\sum_n (1-|\lambda_n|^2)|f(\lambda_n)|^2\simeq
\sum_n 2^n=+\infty$, and $f$ cannot be
in $H^2$ neither in $K^2_I$.

Still $f$ can be written $f=I\overline{\psi}$ with $\psi\in N^+$
(the Smirnov class) and $\psi(0)=0$ (more precisely $\psi(z)
=zf(z)$). It remains to choose $g$ suitably so
that $|f|^2$ is integrable against $|g|^2$.

For this construction we consider the argument of $I$:
\beqa
 \vp(t):=\arg I(e^{it})=\Im \frac{1+e^{it}}{1-e^{it}}
 =\frac{1}{i}\frac{e^{-it}-e^{it}}{|e^{it}-1|^2}
 =-\frac{\cos t}{1-\cos t}
\eeqa 
Observe that $\vp'(t)=\sin t/(1-\cos t)^2$ so that
$\vp$ is strictly increasing on $(0,\pi)$.
Now we consider the intervalls $M_k=[2^{-(k+1)},2^{-k})$. We check
that on these intervals the function $\vp$ increases more than
$2\pi$ ($k$ sufficiently big). 
For this let $t\in (0,\pi/2)$, then there exists
$\xi_t\in (t,2t)$ such that
\beqa
 |\vp(2t)-\vp(t)|=t \frac{\sin\xi_t}{(1-\cos\xi_t)^2}
 \simeq t\frac{\xi_t}{\xi_t^4/4}\simeq \frac{1}{t^2}
\eeqa
(note that $t\le \xi_t\le 2t$).
Since the last expression tends to infinity when 
$t\to 0^+$ there exists an $N\in \N$ such that for every
$n\ge N$, $\vp(2^{-k})-\vp(2^{-(k+1)})\ge 2\pi$.
Hence, by the intermediate value theorem, 
for every $n\ge N$, there exists $t_k\in M_k$ such 
that $\vp(t_k)=0[2\pi]$ and hence $I(e^{it_k})=1$ and $f(e^{it_k})=0$.
Since $t\lmto I(e^{it})$ depends continuously on $t$ outside $0$,
there exists
$\delta_k>0$ such that $|f(e^{it})|\le 1$ for
$t\in F_k:=[t_k-\delta_k,t_k+\delta_k]$.
We will set $E_k=\{e^{it}:t\in F_k\}$
and $\Gamma_k^+:=\{e^{it}:t\in[t_k+\delta_k,t_{k-1}-\delta_{k-1})\}
\subset M_k\cup M_{k-1}$ (we can suppose that $\delta_k$ is
sufficiently small so that $\Gamma_k^+$ is non void).
We will also use the symmetric arc (with
respect to the real axis): $\Gamma_k^-:=\overline{\Gamma_k^+}$ and
$\Gamma:=\Gamma_k^+\cup\Gamma_k^-$.
Then $|\Gamma_k^{\pm}|\lesssim 1/2^k$.
Now for $s<2$ pick $r\in (s,2)$
and let $\alpha_k=1/2^{k/r}$. Define $g$
to be the outer function in $H^2$ such that
\beqa
 |g|=\omega:=\sum_{k}\alpha_k\chi_k+\chi_{\T\setminus\Gamma} \quad
 \text{a.e. }\T,
\eeqa
where $\chi_k=\chi_{\Gamma_k}$ is the characteristic function of
the set $\Gamma_k$ and $\Gamma=\bigcup_n\Gamma_n$.
The function $\omega$ is bounded and $\log$-integrable:
\beqa
 \int_{\T}|\log w| dm=\sum_k|\Gamma_k| |\log \alpha_k|
 \lesssim\sum_k \frac{2}{2^k}\frac{k}{r}\log
 2<\infty.
\eeqa
Hence $g$ is in $\Hi\subset H^2$. We check that $1/g\in H^s$.
\beqa
 \int_{\T}\frac{1}{|g|^s}dm=|\T\setminus\Gamma|
 +\sum |\Gamma_k| {(2^{k/r})^s}
 \lesssim 1+2\sum 2^{k(s/r-1)}
\eeqa
which converges since $s<r$.
Now
\beqa
 \int_{\T}|f|^2|g|^2&=&\int_{\T\setminus\Gamma} \left|
 \frac{1-I}{1-z}\right|^2dm+\sum \alpha_k^2\int_{\Gamma_k}
 \left| \frac{1-I}{1-z}\right|^2dm\\
 &\lesssim& 1+\sum  \frac{8}{(2^{k/r})^2}\frac{1}{2^k}2^{2k}\\
 &\lesssim& 1+8\sum 2^{k(2-1-2/r)}
\eeqa
which converges since $r<2$
(note that for the estimate on $\T\setminus\Gamma$, we have
used the fact that on $E_k$ the function $f$ is bounded by $1$,
and on the remaining closed
arc joining $t_1+\delta_1$ to $2\pi-(t_1+\delta_t)$ it is continuous.
\end{proof}

Observe that the function $g$ constructed in the previous proof
is big (equal to 1) on small sets coming arbitrarily close to 1,
and tending to zero on the remaining sets when $e^{it}\to 1$. In
particular, such a function cannot satisfy the Muckenhoupt 
condition.\\

The examples in the preceding two propositions 
show that it is not always possible to
embed $K^2_I(|g|^2)$ into $K^2_I$ under different conditions
on $|g|$. However in both cases we have the global
integrability condition that appeared 
in Proposition \ref{prop4n},
so that in these cases we can embed $K^p_I(|g|^p)$ into
a $K^{r}_I$ for a suitable $r>1$.

\end{document}